\documentclass[12pt]{article}
\usepackage{amsmath}

\newtheorem{cor}{Corollary}
\newtheorem{theorem}[cor]{Theorem}
\newtheorem{prop}[cor]{Proposition}
\newtheorem{lem}[cor]{Lemma}
\newtheorem{re}[cor]{Remark}
\newtheorem{ex}[cor]{Example}

\newenvironment{remark} {\begin{re} \upshape} {\end{re}}

\newenvironment{proof}{\noindent \emph{Proof.}}  {\vspace{.2cm}}

\begin{document}

\title{Bounds for DNA codes with constant GC-content}

\author{Oliver D. King}
\date{}
\maketitle

\footnotetext{Department of Biological Chemistry and Molecular
Pharmacology, Harvard Medical School,  250 Longwood Avenue,
SGMB-322 Boston, Massachusetts 02115,
\texttt{oliver\_king@hms.harvard.edu}; supported in part by a
fellowship from NIH/NHGRI.}

\maketitle

\begin{abstract}
We derive theoretical upper and lower bounds on the maximum size
of DNA codes of length $n$ with constant GC-content $w$ and
minimum Hamming distance $d$, both with and without the additional
constraint that the minimum Hamming distance between any codeword
and the reverse-complement of any codeword be at least $d$. We
also explicitly construct codes that are larger than the best
previously-published codes for many choices of the parameters $n$,
$d$ and $w$.
\end{abstract}

\section*{Introduction}

Libraries of DNA words satisfying certain combinatorial
constraints have applications to DNA barcoding and DNA computing
(see e.g.\ \cite{cond} and the references therein). The goal is to
design libraries that are as large as possible given the
constraints.

We first review some terminology and notation --- see
\cite{macw,cond} for more context. Let $Z_q$ denote the
$q$-character alphabet $\{0, \ldots, q-1\}$. By a \emph{$q$-ary
word} of length $n$ we mean an element $\mathbf{x}$ of $Z_q^n$,
which we write as $\mathbf{x} = x_1 \cdots x_n$. A \emph{$q$-ary
code} of length $n$ is just a subset of $Z_q^n$, and the elements
of the code are called \emph{codewords}. The \emph{Hamming
distance} $H(\mathbf{x},\mathbf{y})$ between two $q$-ary words
$\mathbf{x}$ and $\mathbf{y}$ of length $n$ is defined to be the
number of coordinates in which they differ, and the \emph{Hamming
weight} of $\mathbf{x}$ is the number of coordinates in which it
is nonzero. The maximum cardinality of a $q$-ary code of length
$n$ for which the minimum Hamming distance between two distinct
codewords is at least $d$ is denoted $A_q(n,d)$. If we also
require each codeword to have Hamming weight $w$ (i.e., that the
code be a \emph{constant-weight code}), the maximum cardinality is
denoted $A_q(n,d,w)$.

A \emph{DNA code} is a $q$-ary code with $q=4$; we identify the
elements 0, 1, 2, 3 $\in Z_4$ with the nucleotides $A$, $C$, $G$,
$T$ (in that order). The \emph{reverse complement} of a DNA word
$\mathbf{x} = x_1 \cdots x_n$ is denoted by $\mathbf{x}^{RC}$, and
is defined to be the word $\overline{x_n} \cdots \overline{x_1}$
where $\overline{x_i}$ is the Watson-Crick complement of $x_i$
(i.e., $\overline{A} = T$, $\overline{T} = A$, $\overline{C}=G$,
and $\overline{G} = C$). By requiring the minimum Hamming distance
between two DNA codewords to be sufficiently large, one can make
it unlikely that a codeword hybridizes to the reverse-complement
of any other codeword. By requiring the minimum Hamming distance
between a DNA codeword and the reverse-complement of a DNA
codeword to be sufficiently large, one can make it unlikely that a
codeword hybridizes to any other codeword or to itself
\cite{cond6}. We denote by $A^{RC}_4(n,d)$ the maximum size of a
DNA code of length $n$ in which $H(\mathbf{x},\mathbf{y}) \geq d$
for all distinct codewords $\mathbf{x}$ and $\mathbf{y}$ and
$H(\mathbf{x},\mathbf{y}^{RC}) \geq d$ for all (not-necessarily
distinct) codewords $\mathbf{x}$ and $\mathbf{y}$. If we also
require each codeword to have Hamming weight $w$ the maximum
cardinality is denoted $A^{RC}_4(n,d,w)$.

The \emph{GC-content} of a DNA word is defined to be the number of
positions in which the word has coordinate $C$ or $G$. It may be
desirable that all codewords in a DNA code have roughly the same
GC-content, so that they have similar melting temperatures (see
e.g.\ \cite{cond6});  $A_4^{GC}(n,d,w)$ and $A_4^{GC,RC}(n,d,w)$
are defined analogously to $A_4(n,d,w)$ and $A_4^{RC}(n,d,w)$,
except that in the former two cases it is the GC-content (rather
than the Hamming weight) of each codeword that is required to be
$w$.

Theoretical upper and lower bounds on $A_4^{RC}(n,d,w)$, with no
restriction on GC-content, are given in \cite{cond}. Explicit
constructions using stochastic local search \cite{cond2,tulpan2}
and a ``template-map'' strategy \cite{cond3} provide lower bounds
on $A_4^{GC}(n,d,w)$ and $A_4^{GC,RC}(n,d,w)$ for a limited range
of parameters $n$, $d$ and $w$. In this paper we derive
theoretical upper and lower bounds on $A_4^{GC}(n,d,w)$ and
$A_4^{GC,RC}(n,d,w)$ for all parameters, and we use lexicographic
constructions to find explicit codes that improve on many of the
lower bounds in \cite{cond3,cond2,tulpan2}.

\section*{Upper bounds}

Before giving upper bounds on the sizes of DNA codes with constant
GC-content, we note some simple special cases:

\begin{prop} For $n>0$, with $0 \leq d \leq n$ and $0 \leq w \leq n$,
\begin{eqnarray}
A_4^{GC}(n,d,0) & = & A_2(n,d) \label{a}\\
A_4^{GC}(n,d,w) & = & A_4^{GC}(n,d,n-w) \label{b}\\
A_4^{GC}(n,n,w) & = & \left\{ \begin{array}{l} 4 \mbox{ if } w = n/2 \\
  3 \mbox{ if } n/3 \leq w < n/2 \mbox{ or } n/2 < w \leq 2n/3 \\
2 \mbox{ if } w < n/3 \mbox{ or } w > 2n/3 \\
\end{array} \right.  \label{c}\\
A_4^{GC,RC}(n,n,w) & = & \left\{ \begin{array}{l} 2 \mbox{ if } w = n/2 \\
  1 \mbox{ if } w \neq n/2 \\
\end{array} \right. \label{d}\\
A_4^{GC}(n,1,w) & = & \binom{n}{w} 2^n \label{e}\\
A_4^{GC,RC}(n,1,w) & = & \left\{
\begin{array}{l} \frac{1}{2}(\binom{n}{w} 2^n -
\binom{n/2}{w/2} 2^{n/2}) \mbox{ if $n$ is even and $w$ is even,} \\
 \frac{1}{2} \binom{n}{w} 2^n \mbox{ if n
is odd or w is odd.}\end{array}\right. \label{f}
\end{eqnarray}
\end{prop}

\begin{proof}
(\ref{a}): Changing all 0's in a binary code to $A$'s and all 1's
to $T$'s gives a Hamming-distance-preserving bijection between the
set of all binary codes of length $n$ and the set of all DNA codes
of length $n$ with constant GC-content 0.

(\ref{b}): Interchange $A$'s with $C$'s, and $T$'s with $G$'s.

(\ref{c}): By (\ref{b}) we may assume $w \leq n/2$. If no two
codewords agree in any position, then there can be at most four
codewords by the pigeonhole principle. Hence $A(n,n,w) \leq 4$ for
all $w$. If there are four codewords none of which agree in any
position, then each of the four nucleotides must occur exactly
once in each of the $n$ positions, so the average GC-content of
the four words is exactly $n/2$. This implies that $A(n,n,w) \leq
3$ for $w < n/2$, since in a code with constant GC-content $w$,
the average GC-content is $w$. If three words each have GC-content
$w < n/3$, then there is some position $j$ in which none of the
words has a $C$ or $G$, and at least two of the three words must
agree in this position (both $A$ or both $T$). Hence $A(n,n,w)
\leq 2$ if $w < n/3$.  The following constructions demonstrate the
reverse inequalities:  For $w=n/2$, the four words $A^wC^w$,
$C^wA^w$, $T^wG^w$ and $G^wT^w$ have pairwise distance $n$; for
$n/3 \leq
 w < n/2$ the three words $C^w A^{n-w}$, $T^{n-w}C^w$ and
$A^{\lfloor (n-w)/2 \rfloor} G^w T^{\lceil (n-w)/2 \rceil}$ have
pairwise distance $n$; for $w<n/3$ the two words $C^{w}A^{n-w}$
and $G^{w}T^{n-w}$ are distance $n$ apart.

(\ref{d}): For $w = n/2$, the two words $A^wC^w$ and $C^wA^w$
satisfy the distance and reverse-complement constraints. For $w
\neq n/2$, the word $C^wA^{n-w}$ satisfies the constraints. These
are the largest sets possible, by (\ref{c}) together with
Theorem 7.

(\ref{e}): This is the total number of DNA words of length $n$ and
GC-content $w$.

(\ref{f}): When $n$ and $w$ are even, there are $\binom{n/2}{w/2}
2^{n/2}$ words with GC-content $w$ that are their own reverse
complements, otherwise there are none.
\end{proof}

\subsection*{Johnson-type bounds}

A code of length $n$ can be \emph{shortened} to a (usually
smaller) code of length $n-1$ without increasing the minimum
Hamming distance, by choosing any character $b \in Z_q$ and any
position $i \in \{1, \ldots, n\}$, keeping just those codewords
that have $b$ in their $i$-th position, and then deleting the
$i$-th position from these codewords \cite{macw}. This procedure
is used in proving the following bounds.

\begin{theorem}
For $0 \leq d \leq n$ and $0 < w < n$,
\begin{eqnarray}
A_4^{GC}(n,d,w) & \leq & \lfloor \frac{2n}{w} A_4^{GC}(n-1,d,w-1)
\rfloor  \label{g} \\
A_4^{GC}(n,d,w) & \leq & \lfloor \frac{2n}{n-w} A_4^{GC}(n-1,d,w)
\rfloor. \label{h}
\end{eqnarray}
\end{theorem}

\begin{proof}
(\ref{g}): In any set of $M$ words with length $n$, minimum
Hamming distance at least $d$ and constant GC-content $w$, there
is some position $i$ in which at least $\lceil wM/2n \rceil$
codewords have nucleotide $C$, or some position $i$ in which at
least $\lceil wM/2n \rceil $ codewords have nucleotide $G$
--- otherwise, the average GC-content would be less than $w$.
Keeping just these codewords, and deleting position $i$, gives a
code with length $n-1$, GC-content $w-1$, and minimum Hamming
distance at least $d$. Inequality (\ref{h}) is analogous, based on
the observation that there is some position with at least $\lceil
(n-w)M/2n \rceil$ $A$'s or $\lceil (n-w)M/2n \rceil$ $T$'s.
\end{proof}

\begin{remark}
Upper bounds on $A_4^{GC}(n,d,w)$ are obtained by repeatedly
applying inequalities (\ref{g}) and (\ref{h}), in any order, until
$n=d$, $n=w$ or $w=0$, at which point (\ref{a})--(\ref{c}) may be
used. (Different orders of applying (\ref{g}) and (\ref{h}) may
result in different bounds.) One may continue using (\ref{h}) even
after $w=0$ (or (\ref{g}) even after $n=w$), until $n=d$, but this
amounts to upper-bounding $A_4^{GC}(n,d,0) = A_2(n,d)$ with the
Singleton bound, $2^{n-d+1}$ (see e.g.\ \cite{brouwer}). Tighter
upper bounds for $A_2(n,d)$ are known for many $n$ and $d$
--- see for example \cite{simon}.
\end{remark}

\begin{theorem} Suppose there is a set of $M$ words of length $n$,
constant GC-content $w$, and minimum Hamming distance at least
$d$. Write $wM = nk + r$ with $0 \leq r < n$. Then
\begin{equation} \label{i}
\begin{array}{rcl} M(M-1)\,d & \leq & (n-r)\,(M^2
-\lfloor\frac{k}{2}\rfloor^2 -\lceil\frac{k}{2}\rceil^2
                         -\lfloor\frac{M-k}{2}\rfloor^2 -\lceil \frac{M-k}{2} \rceil^2
                         ) \vspace{.2cm}\\
         & + & r\,(M^2 - \lfloor\frac{k+1}{2}\rfloor^2 -\lceil\frac{k+1}{2}\rceil^2
            -\lfloor\frac{M-k-1}{2}\rfloor^2
            -\lceil\frac{M-k-1}{2}\rceil^2).
\end{array}
\end{equation}
\end{theorem}

\begin{proof} Let $a_i, c_i, g_i$ and $t_i$ denote the number of
occurrences of $A$, $C$, $G$ and $T$ (respectively) in the $i$-th
position of the $M$ codewords. Note that $\sum_{i=1}^n (c_i + g_i)
= wM$. The sum of the Hamming distances over all $M^2$ ordered
pairs of codewords is $D = \sum_{i=1}^n (M^2 - a_i^2 - c_i^2 -
g_i^2 - t_i^2)$.   Subject only to the constraints that $a_i + c_i
+ g_i + t_i = M$ for each $i$ and that $\sum_{i=1}^n (c_i + g_i) =
wM$, the expression $D$ is maximized when $c_i + g_i$ is as close
as possible to $wM/n$ for each $i$, when $a_i$ is as close as
possible to $t_i$ for each $i$, and when $c_i$ is as close as
possible to $g_i$ for each $i$. This is also true when $a_i$,
$c_i$, $g_i$ and $t_i$ are constrained to be integers, as can be
proved using the same type of argument as in \cite{oster1}, for
example. Hence the right-hand-side of (\ref{i}) is an upper bound
for the sum of the $M^2$ pairwise Hamming distances. For the
left-hand-side, note that since the Hamming distance between
distinct codewords is at least $d$, the sum of the Hamming
distances taken over all $M^2$ ordered pairs of codewords is at
least $M(M-1)\,d$.
\end{proof}

If we relax the constraint that the counts $a_i, c_i, g_i$ and
$t_i$ be integers, Theorem 4 simplifies to the following:
\begin{theorem}
If $2dn > w^2 + 4w(n-w) + (n-w)^2$, then
\begin{equation}
 A_4^{GC}(n,d,w) \leq \frac{2dn}{2dn - (w^2 + 4w(n-w) +
 (n-w)^2)}. \end{equation}
\end{theorem}

\begin{remark}
Versions of the bounds in Theorems 2, 4 and 5 for binary
constant-weight codes \cite{john1,john2} are called Johnson
bounds. Johnson bounds have been generalized to $q$-ary
constant-weight codes \cite{lint,etzion} and to $q$-ary
\emph{constant-composition codes} (where the number of occurrences
of each character in each codeword is prescribed) \cite{oster3}.
They can also be generalized to a setting in which the $q$
characters $\{0, \ldots, q-1 \}$ are partitioned into any number
of subsets, with the total number of occurrences from each subset
specified. Constant-weight codes correspond to the partition
$\{0,\ldots,q-1\} = \{0\} \cup \{1, \ldots, q-1\}$, and
constant-composition codes to the partition $\{0,\ldots,q-1\} =
\{0\} \cup \cdots \cup \{q-1\}$. Our bounds for DNA codes with
constant GC-content correspond to the partition $\{0,1,2,3\} =
\{0,3\} \cup \{1,2\}$.
\end{remark}

\subsection*{Halving bound}
Any upper bound for $A_4^{GC}(n,d,w)$ yields an upper bound for
$A_4^{GC,RC}(n,d,w)$ by the following result, an analogue of the
halving bound for DNA codes with unrestricted GC-content in
\cite{cond}. The same proof works here, since the
reverse-complement of a DNA word has the same GC-content as the
word itself.

\begin{theorem} \label{hb} For $0 < d \leq n$ and $0 \leq w \leq n$,
\begin{equation} A_4^{GC,RC}(n,d,w) \leq  \frac{1}{2}
A_4^{GC}(n,d,w). \label{j} \end{equation}
\end{theorem}

\begin{proof}
If $\{\mathbf{x}_i\}_{i=1}^M$ is a set of $M$ codewords with
constant GC-content $w$, minimum Hamming distance at least $d$,
and with $H(\mathbf{x}_i,\mathbf{x}_j^{RC}) \geq d$ for all $1
\leq i,j \leq M$, then $\{\mathbf{x}_i\}_{i=1}^M \cup
\{\mathbf{x}_i^{RC}\}_{i=1}^M$  is a set of words with constant
GC-content $w$ and minimum Hamming distance at least $d$. This set
has cardinality $2M$ provided that $\{\mathbf{x}_i\}_{i=1}^M \cap
\{\mathbf{x}_i^{RC}\}_{i=1}^M = \emptyset,$ which holds for $d>0$.
\end{proof}

\section*{Lower bounds}

\subsection*{Gilbert-type bounds}

If $\mathcal{C}$ is set of words in $Z_q^n$ with the property that
the Hamming distance between any pair of words in $\mathcal{C}$ is
at least $d$, and if $\mathcal{C}$ is maximal in the sense that no
more points from $Z_q^n$ can be added to $\mathcal{C}$ without
violating this distance constraint, then the balls of Hamming
radius $d-1$ around the points in $\mathcal{C}$ cover all of
$Z_q^n$. This is the idea behind the Gilbert bound for $q$-ary
codes (see e.g.\ \cite{pless}), and a similar argument applies to
constant-weight codes (see e.g.\ \cite{sloane1}). Here we give an
analogue for DNA codes with constant GC-content:

\begin{theorem} For $0 \leq d \leq n$ and $0 \leq w \leq n$,
\begin{equation} A_4^{GC}(n,d,w) \geq \frac{\binom{n}{w} 2^n}
{\sum_{r=0}^{d-1} \sum_{i=0}^{\min \{\lfloor r/2 \rfloor, w,
n-w\}} \binom{w}{i} \binom{n-w}{i} \binom{n-2i}{r-2i} 2^{2i}}.
\label{k} \end{equation}
\end{theorem}

\begin{proof}
The numerator gives the total number of words with GC-content $w$.
The denominator gives the number of these words that have distance
at most $d-1$ from any fixed codeword $\mathbf{x}$. (In the
denominator, $\binom{w}{i} \binom{n-w}{i} \binom{n-2i}{r-2i}
2^{2i}$ is the number of words $\mathbf{y}$ with GC-content $w$
for which $H(\mathbf{x},\mathbf{y})$ is exactly $r$, and for which
there are exactly $w-i$ positions $j$ with $x_j$ and $y_j$ both in
$\{C,G\}$.)
\end{proof}

\begin{remark}
Replacing $d-1$ with $\lfloor (d-1)/2 \rfloor$ as the upper index
of the outer summation in the denominator of (\ref{k}) gives an
upper-bound for $A_4^{GC}(n,d,w)$, since the balls of Hamming
radius $\lfloor (d-1)/2 \rfloor$ centered around codewords must be
disjoint. This is an analogue of the sphere-packing bound for
$q$-ary codes --- see e.g.\ \cite{pless}.

\end{remark}

Now define $V(n,w,d) = \#\{\mathbf{x} \in Z_4^n :  \mathbf{x}
\mbox{ has GC-content } w \mbox{ and }
H(\mathbf{x},\mathbf{x}^{RC})=d\}$. Note that since no nucleotide
is its own complement, $V(n,w,d) = 0$ unless $n$ and $d$ have the
same parity (i.e., are both even or are both odd).

\begin{lem} For $n=2m$ and $d=2e$ even,
\begin{equation} V(2m,w,2e) = \sum_{i=\max \{0,w-m,\lceil (w-e)/2
 \rceil \}}^{\lfloor w/2 \rfloor} \binom{m}{i} \binom{m-i}{w-2i}
 \binom{m - w + 2i}{e - w + 2i} 2^{m + 2w - 4i}; \label{l} \end{equation}

For $n=2m+1$ and $r=2e+1$ odd, \begin{equation} V(2m+1,w,2e+1) =
V(2m,w,2e) + V(2m,w-1,2e). \label{m} \end{equation}
\end{lem}

\begin{proof}
In (\ref{l}), the index $i$ ranges over the number of positions $j
\leq m$ for which both $x_j$ and $x_{2m-j+1}$ belong to $\{C,G\}$.
There are $\binom{m}{i}$ ways to select these positions, and
$\binom{m-i}{w-2i} 2^{w-2i}$ ways to select the positions for the
remaining $w-2i$ occurrences of $C$'s or $G$'s. There are then
$m-w+i$ positions $j \leq m$ for which both $x_j$ and $x_{2m-j+1}$
belong to $\{A,T\}$. Note that the $j$-th coordinate of
$\mathbf{x}$ necessarily differs from the $j$-th coordinate of
$\mathbf{x}^{RC}$ in the $w-2i$ positions $j \leq m$ for which one
of $x_j$ and $x_{2m-j+1}$ is in $\{A,T\}$ and the other is in
$\{C,G\}$, so there are $\binom{m - w + 2i}{e - w + 2i}$ ways to
choose the remaining $e-w + 2i$ positions $j \leq m$ in which
$x_j$ differs from the complement of $x_{2m+1-j}$. After all these
choices have been made, there are two choices for the nucleotide
in each position $j \leq m$; for the $m-w+2i$ positions $j \leq m$
for which $x_j$ and $x_{2m-j+1}$ both belong to $\{C,G\}$ or both
belong to $\{A,T\}$, the nucleotide at $x_{2m-j+1}$ is forced by
the choice of $x_j$; for the other $w-2i$ positions $j \leq m$,
there are two choices for the nucleotide $x_{2m-j+1}$.

In (\ref{m}), the first summand gives the number of words with
$x_{m+1} \in \{A,T\}$ and the second summand gives the number of
words with $x_{m+1} \in \{C,G\}$.

\end{proof}

\begin{theorem} \label{glb} For $0 \leq d \leq n$ and $0 \leq w \leq n$,
\begin{equation} A_4^{GC,RC}(n,d,w) \geq \frac{ \sum_{r=d}^{n} V(n,d,r)} {2
\sum_{r=0}^{d-1} \sum_{i=0}^{\min \{\lfloor r/2 \rfloor, w, n-w\}}
\binom{w}{i} \binom{n-w}{i} \binom{n-2i}{r-2i} 2^{2i}}.
\label{n}\end{equation}
\end{theorem}

\begin{proof}
The numerator gives the total number of words with GC-content $w$
that have distance at least $d$ from their reverse-complements,
and the denominator gives an upper-bound on the number of these
words that have distance at most $d-1$ from any fixed codeword.
(The denominator is an upper-bound rather than an exact count,
because the balls of radius $d-1$ around a word and its
reverse-complement might overlap, and because when counting the
number of words in these balls we may be including some words
$\mathbf{y}$ that do not satisfy the condition
$H(\mathbf{y},\mathbf{y}^{RC}) \geq d$.)
\end{proof}

\subsection*{Lexicographic codes}
See \cite{sloane2} for an introduction to lexicographic codes. The
idea is that all words in $Z_q^n$ are listed in lexicographic
order, i.e., with $\mathbf{x} = x_1 \cdots x_n$ listed before
$\mathbf{y} = y_1 \cdots y_n$ if $x_i < y_i$, where $i$ is the
first position in which $\mathbf{x}$ and $\mathbf{y}$ differ.
Then, starting with the empty code, one proceeds down this list
and adds to the code any word whose addition does not violate any
of the combinatorial constraints. (Ordinarily these would be a
Hamming distance and possibly a Hamming weight constraint, but
GC-content and reverse-complement Hamming distance constraints can
be enforced as well.)  Since the resulting lexicographic codes can
accommodate no more codewords without a constraint being violated,
they meet or exceed Gilbert-type lower bounds; they often do much
better \cite{sloane2}. There are many variants of the standard
lexicographic construction, for example the words may be ordered
as a Gray code, or one may start with an arbitrary codeword as a
seed rather than with the empty code \cite{sloane1}. We used three
variants, singly and in combination, to construct DNA codes with
the desired constraints:

(i) We used different orderings of the characters $A$, $C$, $G$
and $T$ when putting the $4^n$ DNA words of length $n$ in
lexicographic order. There are $4! = 24$ orderings of the four
characters, but because of the symmetry between $A$ and $T$ and
between $C$ and $G$, only 6 of these 24 orderings need to be
considered.

(ii) We used offsets, as in \cite{oster1}: one starts at an
arbitrary place in the list of words rather than at the beginning,
and loops back around to the beginning of the list when the end is
reached.

(iii) We used a ``factored'' ordering of the DNA words. The $2^n$
binary words of length $n$ were listed in lexicographic order,
$\mathbf{u}_1 = 0 \cdots 0, \ldots , \mathbf{u}_{2^n} = 1 \cdots
1$. As in \cite{cond}, we define a mapping $\odot$ from pairs of
binary words of length $n$ to DNA words of length $n$, given by
$\mathbf{x} \odot \mathbf{y} = \mathbf{z}$ where $z_i = A$ if
$x_i=0$ and $y_i=1$; $z_i=C$ if $x_i=1$ and $y_i=0$; $z_i = G$ if
$x_i=y_i=1$; and $z_i=T$ if $x_i = y_i = 0$. Note that $\odot$ is
a bijection, and that the Hamming weight of $\mathbf{x}$ is equal
to the GC-content of $\mathbf{z}$.  We ordered the $4^n$ DNA words
so that $\mathbf{u}_i \odot \mathbf{u}_j$ comes before
$\mathbf{u}_k \odot \mathbf{u}_m$ if $i<k$ or if $i=k$ and $j<m$.

When combining variants (ii) and (iii) above, two offsets can be
used: one for the binary words in the first slot of $\mathbf{x}
\odot \mathbf{y}$, and another for those in the second slot.

We used the above three approaches to construct DNA codes with
constant GC-content, both with and without the reverse-complement
constraint, for a variety of parameters $n$, $d$ and $w$.  Using
offsets of zero and an average of about ten random offsets, we
found codes that are larger than the codes given in
\cite{cond3,cond3,tulpan2} for many choices of parameters. The
sizes of the lexicographic codes are given in Tables 1 and 2, and
the offsets used to generate these codes are given in Tables 3 and
4.

\subsection*{Product bounds}

The lexicographic constructions described above do not scale well
to large $n$. One can avoid the burden of explicitly computing
distances between all pairs of codewords (and also the burden of
explicitly listing all codewords) by using modifications of
algebraic constructions such as linear codes. For example, a DNA
code with minimum Hamming distance at least $d$ and constant
GC-content $w$ can be constructed by taking any linear code over
$Z_4$ (or the Galois field $\mathbf{F}_4$ \cite{cald} or the
Kleinian four-group \cite{hohn}) that has minimum Hamming distance
$d$, and selecting only those codewords with exactly $w$
occurrences of two fixed characters.

In this section we give lower bounds for DNA codes that are
constructed from binary codes, binary constant-weight codes, and
ternary constant-weight codes, for which a variety of algebraic
constructions are known (e.g.\ \cite{macw,sloane1,oster1}).

Note that the reverse-complement operator RC can be viewed as the
composition of two (commuting) operators R and C, where R maps
$x_1 \cdots x_n$ to $x_n \cdots x_1$ and C replaces each
coordinate $x_i$ with its complement $\overline{x_i}$. We state
the product bounds below in terms of constraints on $R$ rather
than on $RC$ to make the arguments cleaner. (This approach was
used in \cite{cond}.) The values $A_q^{R}(n,d,w)$ and
$A_4^{GC,R}(n,d,w)$ are defined in the same manner as
$A_q^{RC}(n,d,w)$ and $A_4^{GC,RC}(n,d,w)$, but with the
constraint that $H(\mathbf{x},\mathbf{y}^R) \geq d$ for all
codewords $\mathbf{x}$ and $\mathbf{y}$ in place of the constraint
that $H(\mathbf{x},\mathbf{y}^{RC}) \geq d$. Bounds on
$A_4^{GC,R}(n,d,w)$ can be used to derive bounds for
$A_4^{GC,RC}(n,d,w)$ using the following result:

\begin{prop} For $0 \leq d \leq n$ and $0 \leq w \leq n$,
\begin{equation}
A_4^{GC,RC}(n,d,w) = A_4^{GC,R}(n,d,w) \mbox{ if $n$ is even},
\label{o}
\end{equation}
\begin{equation} A_4^{GC,R}(n,d+1,w)  \leq
A_4^{GC,RC}(n,d,w) \leq A_4^{GC,R}(n,d-1,w) \mbox{ if $n$ is odd}.
\label{p}
\end{equation}
\end{prop}

\begin{proof}
The analogous result for DNA codes with unrestricted GC-content
was proved in \cite{cond}, and essentially the same proof works
here. Given a set of codewords of length $n$, if we replace all
the entries in any subset of the positions by their complements,
the GC-content of each codeword is preserved, as is the Hamming
distance between any pair of codewords. The Hamming distance
between a codeword and the reverse or reverse-complement of
another codeword is not in general preserved, but if $n$ is even
and we replace the first $n/2$ coordinates of each codeword
$\mathbf{x}_i$ by their complements to form a new word
$\mathbf{y}_i$, then $H(\mathbf{x}_i,\mathbf{x}_j^R) =
H(\mathbf{y}_i,\mathbf{y}_j^{RC})$ for all codewords
$\mathbf{x}_i$ and $\mathbf{x}_j$. Similarly, if $n$ is odd and we
replace the first $(n-1)/2$ coordinates of each codeword
$\mathbf{x}_i$ by their complements to form $\mathbf{y}_i$, then
$|H(\mathbf{x}_i,\mathbf{x}_j^R) -
H(\mathbf{y}_i,\mathbf{y}_j^{RC})| \leq 1$.
\end{proof}

\begin{theorem}
For $0 \leq d \leq n$ and $0 \leq w \leq n$,
\begin{eqnarray}
A_4^{GC}(n,d,w) & \geq &  A_2(n,d,w) \cdot A_2(n,d)  \label{q}\\
A_4^{GC,R}(n,d,w) & \geq & A_2^{R}(n,d,w) \cdot A_2(n,d) \label{r} \\
A_4^{GC,R}(n,d,w) & \geq & A_2(n,d,w) \cdot A_2^R(n,d) \label{s} \\
A_4^{GC}(n,d,w) & \geq & A_3(n,d,w) \cdot A_2(n-w,d) \label{t} \\
A_4^{GC,R}(n,d,w) & \geq & A_3^R(n,d,w) \cdot A_2(n-w,d) \label{u}\\
A_4^{GC,R}(n,d,w) & \geq & A_3(n,d,w) \cdot A_2^R(n-w,d) \label{v}
\end{eqnarray}
\end{theorem}

\begin{proof}
 For (\ref{q}) and (\ref{r}), note that if $\mathcal{B}_1$ is a set of binary words with
 length $n$, Hamming weight $w$ and minimum Hamming distance $d$,
 and if $\mathcal{B}_2$ is a set of binary words with length $n$ and
 minimum
 Hamming distance $d$, then $\mathcal{D} = \{\mathbf{x} \odot \mathbf{y} : \mathbf{x} \in
 \mathcal{B}_1 \mbox{ and } \mathbf{y} \in \mathcal{B}_2\}$
 is a set of DNA words with length $n$, GC-content $w$ and minimum
 Hamming distance $d$. If, in addition, $H(\mathbf{x}_1,\mathbf{x}_2^R)
\geq d$ for all $\mathbf{x}_1, \mathbf{x}_2 \in \mathcal{B}_1$,
then $H(\mathbf{z}_1,\mathbf{z}_2^R) \geq d$ for all
$\mathbf{z}_1, \mathbf{z}_2 \in \mathcal{D}$ as well, since
$H(\mathbf{x}_1 \odot \mathbf{y}_1, (\mathbf{x}_2 \odot
\mathbf{y}_2)^R) = H(\mathbf{x}_1 \odot \mathbf{y}_1,
\mathbf{x}_2^R \odot \mathbf{y}_2^R) \geq H(\mathbf{x}_1,
\mathbf{x}_2^R) \geq d$. Inequality (\ref{s}) is proved in the
same manner as (\ref{r}).

For (\ref{t})--(\ref{v}) we first define a function $\oslash$ that
maps a pair consisting of ternary word $\mathbf{x}$ of length $n$
and Hamming weight $w$, and a binary word $\mathbf{y}$ of length
$n-w$, to a DNA word $\mathbf{z} = \mathbf{x} \oslash \mathbf{y}$
of length $n$. This map is defined by $z_i = C$ if $x_i = 1$; $z_i
= G$ if $x_i = 2$; $z_i=A$ if $x_i$ is the $j$-th zero-entry in
$\mathbf{x}$ and $y_j=0$; and $z_i=T$ if $x_i$ is the $j$-th
zero-entry in $\mathbf{x}$ and $y_j=1$. The argument now proceeds
as for (\ref{q})--(\ref{s}).
\end{proof}

\begin{remark}
Lower bounds for $A_2(n,d,w)$ can be found in \cite{sloane1},
lower bounds for $A_2(n,d)$ in \cite{brouwer,simon}, and lower
bounds for $A_3(n,d,w)$ in \cite{oster1}. The bounds on ternary
constant-weight codes in \cite{oster1} also apply directly to DNA
codes with constant C-content over the three-letter alphabet
$\{A,C,T\}.$ This restricted alphabet is used by some researchers
to reduce the probability of individual codewords having
``secondary structure'' such as hairpin loops \cite{mir,faul} ---
note also that if $\mathbf{x}$ and $\mathbf{y}$ are DNA words over
$\{A,C,T\}$ with C-content at least $d$, the reverse-complement
Hamming distance constraint $H(\mathbf{x},\mathbf{y}^{RC}) \geq d$
is automatically satisfied.
\end{remark}

\begin{remark}
Inequalities (\ref{q})--(\ref{s}) are analogues of the product
bounds for DNA codes with unrestricted GC-content in \cite{cond};
(\ref{q}) is also a generalization of the ``template-map''
construction used in \cite{cond3} for codes with constant
GC-content --- in that construction, a constant-weight binary code
acts as the ``template'' (corresponding to the first factor in
(\ref{q})), and the same constant-weight binary code, with at most
two words of other weights added in, acts as the ``map''
(corresponding to the second factor in (\ref{q})). This gives a
DNA code of size no larger than $A_2(n,d,w) \cdot A_2(n,d)$, and
when $A_2(n,d,w) + 2 < A_2(n,d)$ this gives a strictly smaller
code (e.g.,  $A_2(n,2,w) = \binom{n}{w}$, which can be much less
than $A_2(n,2) = 2^{n-1}$). But for the parameters $w = d \approx
n/2$ considered in \cite{cond3}, this difference can be
inconsequential; in particular, $A_2(n,n/2,n/2) = A_2(n,n/2) - 2 =
2n -2$ whenever a Hadamard matrix of order $n$ exists
\cite{semakov}, i.e.\ for all $n$ divisible by 4 up to at least $n
= 424$. Note that even when optimal binary codes are used as
factors, the lower bounds derived from product codes are not in
general tight --- for instance,
$A_2(12,6,6) \cdot A_2(12,6) = 22 \cdot 24 = 528$,
while we constructed a lexicographic code showing that
$A_4^{GC}(12,6,6) \geq 736$. In fact, product codes do not even
meet the Gilbert-type lower bound for $A_4^{GC}(2w,w,w)$ when $w$
is sufficiently large: replacing the denominator in (\ref{k}) with
the upper-bound $w \binom{2w}{w-1} 3^{w-1}$ for the number of
words with Hamming distance at most $w-1$ from a fixed codeword gives
$A_4^{GC}(2w,w,w) \geq 3(4/3)^{w}(w+1)/w^2$; the product-code
construction gives a code of size at most $A_2(2w,w,w) \cdot
A_2(2w,w) \leq (4w-2)4w$. (The ``template-code'' construction used
in \cite{arita,arita2} is similar to the template-map construction
discussed above, but with an additional constraint to prevent
codewords from hybridizing to concatenations of other codewords.)

\end{remark}

Below we show that product codes can be optimal when $d=2$:

\begin{theorem}
For $0 \leq w \leq n$,
\begin{equation}
A_4^{GC}(n,2,w) = \binom{n}{w} 2^{n-1}. \label{w}
\end{equation}
\end{theorem}

\begin{proof}
In one direction we have $A_4^{GC}(n,2,w)  \geq A_2(n,2,w) \cdot
A_2(n,2)$ by (\ref{q}). Note that $A_2(n,2,w) = \binom{n}{w}$
since the Hamming distance between two distinct binary words of
the same weight is at least two; note also that $A_2(n,2) =
2^{n-1}$, since the first $n-1$ coordinates can be arbitrary with
the last coordinate used as a parity check bit (see e.g.\
\cite{pless}).

In the other direction, $A_4^{GC}(w,2,w) = A_4^{GC}(w,2,0) =
A_2(w,2) = 2^{w-1} = \binom{w}{w} 2^{w-1}$, and if
$A_4^{GC}(n,2,w) \leq \binom{n}{w} 2^{n-1}$ for some $n \geq w$
then by (\ref{h}) we have $A_4^{GC}(n+1,2,w) \leq 2(n+1-w)/(n+1)
\binom{n}{w} 2^{n-1} = \binom{n+1}{w} 2^n$. Hence by induction
$A_4^{GC}(n,2,w) \leq \binom{n}{w} 2^{n-1}$ for all $n \geq w$.
\end{proof}

\begin{theorem}
For $0 \leq w \leq n$ and $n$ even,
\begin{equation} A_4^{GC,RC}(n,2,w) = \binom{n}{w} 2^{n-2}.
\label{x}
\end{equation}
\end{theorem}

\begin{proof}
By (\ref{k}), $A_4^{GC,RC}(n,2,w) \leq \frac{1}{2} A_4^{GC}(n,2,w)
= \frac{1}{2} \binom{n}{w} 2^{n-1} = \binom{n}{w} 2^{n-2}$. For
$n$ even, $A_4^{GC,RC}(n,2,w) = A_4^{GC,R}(n,2,w)$ by (\ref{o}),
and $A_2^R(n,2) = 2^{n-2}$ by Theorem 4.5 of \cite{cond}. Thus by
the product bound $A_4^{GC,R}(n,d,w) \geq A_2(n,2,w) \cdot
A_2^R(n,2) = \binom{n}{w} 2^{n-2}$. (Here is an alternate argument
showing $A_2^R(n,2) = 2^{n-2}$ for $n$ even: when $n$ is even, the
set of all $2^{n-1}$ binary words of odd Hamming weight contains
no palindromes, and the reverse of a binary word of odd weight has
odd weight, so these $2^{n-1}$ words break up into $2^{n-2}$ pairs
$\{\mathbf{x},\mathbf{x}^R\}$; taking one word from each pair
shows that $A_2^R(n,2) \geq 2^{n-2}$, since the Hamming distance
between two distinct binary words of odd weight is at least two;
equality follows from a halving bound, $A_2^R(n,2) \leq
\frac{1}{2}A_2(n,2) = 2^{n-2}$ \cite{cond}.
\end{proof}

\section*{Tables}

Lower bounds for $A_4^{GC,RC}(n,d,w)$, derived from codes
constructed using stochastic local search, are given in
\cite{cond2} and \cite{tulpan2} for $n \leq 12$ ($n$ even) with $d
\leq n$ and $w = n/2$. In Tables 1 and 2 we give lower bounds for
$A_4^{GC,RC}(n,d,w)$ and $A_4^{GC}(n,d,w)$ derived from
lexicographic constructions for these same parameters.  Our bounds
are at least as large as those in \cite{cond3,cond2,tulpan2} for
all parameters except the five cases marked with asterisks;
those that are strictly larger (or for which no bounds were given)
are underlined. (Our bound on $A_4^{GC}(n,d,w)$ is not underlined
if it is equal to twice the bound on $A_4^{GC,RC}(n,d,w)$ given in
 \cite{cond3,cond2,tulpan2}, since the former bound is then
implied by the latter using the halving bound.) Entries followed
by periods are optimal, as the lower bounds are equal to the upper
bounds computed using Theorems 2, 4 and 7 (the Johnson-type bounds
and the halving bound). \vspace{12pt}

\noindent Guide to superscripts in Tables 1 and 2:

\noindent $a$. Not explicitly constructed lexicographically; value from
Theorem 16.

\noindent $b$. Not explicitly constructed lexicographically; value
from Theorem 15.

\noindent $*$. Larger code constructed using stochastic local
search in \cite{tulpan2} (size given in superscript).

\vspace{.4cm}

\begin{center}
Table 1. Lower bounds for $A_4^{GC,RC}(n,d,w)$ with $n \leq 12$
($n$ even), $d \leq n$ and $w= n/2$. \vspace{5pt}
\setlength{\tabcolsep}{6pt}
\begin{tabular}{c|ccccccccccc}
$n \backslash d$ & 2 & 3 & 4 & 5 & 6 & 7 & 8 & 9 & 10 & 11 & 12 \\
\hline
4   & 24. & 6. & 2. & - & - & - & - & - & - & - & - \\
6   & \underline{320}. & 39$^{*41}$ & \underline{16} & 4. & 2. & - & - & - & - & - & - \\
8  & \underline{4480}. & 384$^{*390}$ & 112 & 25$^{*26}$ & 10$^{*12}$ & 2. & 2. & - & - & - & - \\
10 & \underline{64512}. & \underline{4084} & \underline{795} & \underline{166} & \underline{46} & 15 & 6 & 2. & 2. & - & - \\
12 & \underline{946176}$.^a$ & \underline{49764} & \underline{8704} & \underline{1362} & \underline{306} & \underline{81} & \underline{27} & \underline{10} & 4. & 2. & 2. \\
 \end{tabular}
 \end{center}
 \vspace{.4cm}

 \begin{center}
 Table 2. Lower bounds for $A_4^{GC}(n,d,w)$ with $n \leq 12$ ($n$
 even), $d \leq n$ and $w = n/2$.
 \vspace{5pt} \setlength{\tabcolsep}{5pt}
 \begin{tabular}{c|ccccccccccc}
$n \backslash d$ & 2 & 3 & 4 & 5 & 6 & 7 & 8 & 9 & 10 & 11 & 12 \\
\hline
4  & 48. & 12. & 4. & - & - & - & - & - & - & - & - \\
6  & \underline{640}. & \underline{96} & \underline{40}. & 8 & 4. & - & - & - & - & - & - \\
8  & \underline{8960}. & \underline{832} & 224 & \underline{56} & 20$^{*24}$ & \underline{5}. & 4. & - & - & - & - \\
10 & \underline{129024}. &  \underline{9344} & \underline{1676} & \underline{360} & \underline{96} & \underline{32} & \underline{16.} & \underline{5}. & 4. & - & - \\
12 & \underline{1892352}$.^b$ & \underline{112640} & \underline{17408} & \underline{2992} & \underline{736} & \underline{177} & \underline{68} & \underline{22} & 8 & 4. & 4. \\
\end{tabular}
\end{center}

\vspace{.4cm}

\begin{remark} In \cite{cond2} and \cite{cond3}, lower bounds for
$A_4^{GC}(n,d,w)$ are given for $4 \leq n \leq 20$ ($n$ odd or
even) with $w = d = \lfloor n/2 \rfloor.$ Though not covered in
Table 2, we also improved upon these bounds for $n=5$, 7, 9, 11
and 13--20 using lexicographic constructions.
\end{remark}

Tables 3 and 4 record the nucleotide-orderings and the offsets
used in constructing the lexicographic codes whose sizes are given
in Tables 1 and 2. Entries are written either in the form
$\mbox{offset}_1 \odot \mbox{offset}_2$ for the ``factored''
lexicographic variant, or as $\mbox{offset}^{k}$, with the
superscript $k \in \{1,2,3,4,5,6\}$ indicating the
nucleotide-ordering used, as follows:  $1 = (A < C < G < T);\,\, 2
= (C < G < A < T);\,\, 3 = (A < T < C < G);\,\, 4 = (C < A < T <
G);\,\, 5= (C < A < G < T);\,\, 6 = (A < C < T < G)$. Note that we
list all offsets in base-16 rather than base-4 or base-2 for
compactness, and that the offset need not itself be a codeword
since it may not satisfy the GC-content constraint or it may be
too close to its own reverse-complement. (We re-used seeds for our
random-number generator, which is why some of the same ``random''
offsets appear for more than one entry.)

\vspace{.4cm}

\begin{center}
 Table 3. Offsets used to generate lexicographic codes giving lower
bounds in Table 1. \footnotesize \vspace{5pt}
\setlength{\arraycolsep}{5pt}
 $\begin{array}{c|ccccccccccc}
n \backslash d & 2 & 3 & 4 & 5 & 6 & 7 & 8 & 9 & 10 & 11 & 12 \\
\hline
4  &  59^1    & 59^2   & 0^1 & - & -  & -   & -  & -  & - & - & - \\
6  &  0^1     & \mbox{42d}^4   & 12 \odot 19 &   \mbox{bfc}^2   & 0^1  & -  & -  & - & - & - & - \\
8  & 5021^1  & \mbox{44dd}^2  & \mbox{4e} \odot 95 &  \mbox{d3de}^5   & \mbox{90a5}^5  & 0^1 & 0^1 & - & - & - & - \\
10 & 0^1  & 0^5 & \mbox{bfc99}^1 &  0^5  & 0^1  & \mbox{c0d96}^1 & \mbox{c54c6}^2  & 0^1 & 0^1 & - & - \\
12 & - & 0 \odot 0  & 0 \odot 0 & 0^2 & \mbox{4121c8}^4 & 0^5 & 0^2 & \mbox{96c697}^1 & \mbox{96c697}^1 & 0^1 & 0^1  \\
 \end{array}$
 \end{center}

 \vspace{.4cm}

 \begin{center}
 Table 4. Offsets used to generate lexicographic codes giving lower
 bounds in Table 2. \footnotesize
\vspace{5pt} \setlength{\arraycolsep}{6pt}
 $\begin{array}{c|ccccccccccc}
n \backslash d & 2 & 3 & 4 & 5 & 6 & 7 & 8 & 9 & 10 & 11 & 12 \\
\hline
4  & 0^1 & 0^1 & 0^1 & - & - & - & - & - & - & - & - \\
6  & 0^1 & 0^2 & 434^1 & 0^1 & 0^1 & - & - & - & - & - & - \\
8  & 0^1 & 5021^2 & 0 \odot 0 & \mbox{2d} \odot 23 & \mbox{90f6}^1 & 0^1 & 0^1 & - & - & - & - \\
10 & 0^1 & 0 \odot 0 & 0^2 & 0 \odot 0 & 0 \odot 0 & 0 \odot 0 & \mbox{c8e60}^5 & \mbox{3792d}^2 & 0^1 & - & - \\
12 & - & 0 \odot 0 & 0 \odot 0 & 0 \odot 0 & 0 \odot 0 & \mbox{c8e605}^1 & 994 \odot \mbox{70b} & 0 \odot 0 & 0^2 & 0^1 & 0^1 \\
\end{array}$
\end{center}


\end{document}